\documentclass[11pt]{article}
\usepackage{mathrsfs,amsmath,amssymb,amsfonts}
\numberwithin{equation}{section}
\begin{document}
\title{Direct and inverse theorems for Bernstein operators with inner
singularities}
\author{Wen-ming Lu$^1$ and Lin Zhang$^2$\thanks{Corresponding author. \textsl{E-mail address}:
\texttt{linyz@zju.edu.cn}(L.-Zhang);\texttt{lu\_wenming@163.com}(W.-Lu).}}
\date{\it $^1$School of Science, Hangzhou Dianzi University,
Hangzhou, 310018 P.R. China\\$^2$Department of Mathematics, Zhejiang
University, Hangzhou, 310027 P.R. China} \maketitle
\mbox{}\hrule\mbox{}\\[0.5cm]
\textbf{Abstract}\\[-0.2cm]^^L
We introduce a new type of Bernstein operators, which can be used to
approximate the functions with inner singularities. The direct and
inverse results of the weighted approximation of this new type of
combinations are obtained.\\~\\
%---------------------------------------------------------------------------------------------------------
\mbox{}\hrule\mbox{}
%------------------------------------------------------------------------------------------
\section{Introduction}
The set of all continuous functions, defined on the interval $I$, is
denoted by $C(I)$. For any $f\in C([0,1])$, the corresponding
Bernstein operators are defined as follows:
$$B_n(f,x):=\sum_{k=0}^nf(\frac{k}{n})p_{n,k}(x),$$
where
$$p_{n,k}(x)={n \choose k}x^k(1-x)^{n-k}, \ k=0,1,2,\ldots,n, \ x\in[0,1].$$
Recently Felten showed the following two theorems in
\cite{Felten}:\\~\\
\textbf{Theorem A.} Let $\varphi(x)=\sqrt{x(1-x)}$ and let
$\phi:[0,1] \longrightarrow R,\ \phi \neq 0$ be an admissible
step-weight function of the Ditzian-Totik modulus of
smoothness(\cite{Totik}) such that $\phi^2$ and $\varphi^2/\phi^2$
are concave. Then, for $f \in C[0,1]$ and $0< \alpha <2,\
|B_n(f,x)-f(x)| \leqslant
\omega_\phi^{2}(f,n^{-1/2}{\frac {\varphi(x)}{\phi(x)}}).$\\~\\
\textbf{Theorem B.} Let $\varphi(x)=\sqrt{x(1-x)}$ and let
$\phi:[0,1] \longrightarrow R,\ \phi \neq 0$ be an admissible
step-weight function of the Ditzian-Totik modulus of smoothness such
that $\phi^2$ and $\varphi^2/\phi^2$ are concave. Then, for $f \in
C[0,1]$ and $0< \alpha <2,\ |B_n(f,x)-f(x)|=O((n^{-1/2}{\frac
{\varphi(x)}{\phi(x)}})^\alpha)$ implies
$\omega_\phi^{2}(f,t)=O(t^\alpha).$
\\~\\
Approximation properties of Bernstein operators have been studied
very well (see \cite{Berens}, \cite{Della},
\cite{Totik}-\cite{Lorentz}, for example).
In order to approximate the functions with singularities, Della
Vecchia et al. \cite{Della} introduced some
kinds of modified Bernstein operators. Throughout the paper, $C$
denotes a positive constant independent of $n$ and $x$, which may be
different in different cases.\\
Let $\phi:[0,1] \longrightarrow R,\ \phi \neq 0$ be an admissible
step-weight function of the Ditzian-Totik modulus of smoothness,
that is, $\phi$ satisfies the following conditions:\\
(\textbf{I}) For every proper subinterval $[a,b] \subseteq [0,1]$
there exists a constant $C_1 \equiv C(a,b)>0$ such that $C_{1}^{-1}
\leqslant
\phi(x) \leqslant C_1$ for $x \in [a,b].$\\
(\textbf{II}) There are two numbers $\beta(0) \geqslant 0$ and
$\beta(1) \geqslant 0$ for which
\begin{align*}
\phi(x)\thicksim \left\{
\begin{array}{lrr}
x^{\beta(0)},&&\ \mbox{as}\ x \rightarrow 0+,    \\
(1-x)^{\beta(1)},&&\ \mbox{as}\ x\rightarrow 1-.
              \end{array}
\right.\\
\end{align*}
($X \sim Y$ means $C^{-1}Y \leqslant X \leqslant CY \mbox{\ for\
some}\ C$).\\
Combining conditions (\textbf{I}) and (\textbf{II}) on $\phi,$ we can deduce that\\
\begin{align*}
C^{-1}\phi_2(x) \leqslant \phi(x) \leqslant C\phi_2(x),\ x \in
[0,1],
\end{align*}
where $\phi_2(x) = x^{\beta(0)}(1-x)^{\beta(1)}.$

%---------------------------------------------------------------------------------------------------------
\section{The Main Results}
Let ${\bar{w}}(x)=|x-\xi|^\alpha,\ 0<\xi<1,\ \alpha>0$ and
$C_{\bar{w}}:= \{{f \in C([0,1] \setminus \{\xi\})
:\lim\limits_{x\longrightarrow\xi}(\bar{w}f)(x)=0 }\}$. The norm in
$C_{\bar{w}}$ is defined by
$\|f\|_{C_{\bar{w}}}:=\|{\bar{w}}f\|=\sup\limits_{0\leqslant
x\leqslant 1}|({\bar{w}}f)(x)|$. Define
\begin{eqnarray*}
W_\phi^{2}:= \{f \in C_{\bar{w}}:f' \in A.C.((0,1)),\
\|{\bar{w}}\phi^2f''\|<\infty\},\\
W_{\bar{w},\lambda}^{2}:= \{f \in
C_{\bar{w}}:f' \in A.C.((0,1)),\
\|{\bar{w}}\varphi^{2\lambda}f''\|<\infty\}.
\end{eqnarray*}
For $f \in C_{\bar{w}}$,
 the weighted modulus of smoothness is defined
by
\begin{align*}
\omega_\phi^{2}(f,t)_{\bar{w}}:=\sup_{0<h\leqslant t}\sup_{0
\leqslant x \leqslant 1}|{\bar{w}}(x)\triangle_{h\phi(x)}^{2}f(x)|,
\end{align*}
where
\begin{eqnarray*}
\Delta_{h\phi}^{2}f(x)=f(x+h\phi(x))-2f(x)+f(x-h\phi(x)),
\end{eqnarray*}
and $\varphi(x)=\sqrt{x(1-x)},\ \delta_n(x)=\varphi(x)+{\frac{1}{\sqrt{n}}}.$\\
Let
\begin{eqnarray*}
\psi(x)=\left\{
\begin{array}{lrr}
10x^3-15x^4+6x^5, &&0< x <1, \\
0,   &&x \leqslant0,  \\
1,  &&x \geqslant 1.
             \end{array}
\right.
\end{eqnarray*}
Obviously, $\psi$ is non-decreasing on the real axis, $\psi\in
C^2((-\infty,+\infty)),\ \psi^{(i)}(0)=0,$ $i=0,1,2.\
\psi^{(i)}(1)=0,\ i=1,2$ and $\psi(1)=1.$ Further, let
\begin{eqnarray*}
x_{1}=\frac{[n\xi-2\sqrt{n}]}{n},\ x_{2}=\frac{[n\xi-\sqrt{n}]}{n},\ x_{3}=\frac{[n\xi+\sqrt{n}]}{n},\ x_{4}=\frac{[n\xi+2\sqrt{n}]}{n},
\end{eqnarray*}
and
\begin{eqnarray*}
{\bar{\psi}}_{1}(x)=\psi(\frac{x-x_{1}}{x_{2}-x_{1}}),\ {\bar{\psi}}_{2}(x)=\psi(\frac{x-x_{3}}{x_{4}-x_{3}}).
\end{eqnarray*}
Consider
\begin{eqnarray*}
P(x):={\frac{x-x_{4}}{x_{1}-x_{4}}}f(x_{1}) + {\frac{x_{1}-x}{x_{1}-x_{4}}}f(x_{4}),
\end{eqnarray*}
the linear function joining the points $(x_{1},f(x_{1}))$ and
$(x_{4},f(x_{4}))$. And let
\begin{eqnarray*}
{\bar{F}}_{n}(f,x):={\bar{F}}_{n}(x)=f(x)(1-{\bar{\psi}}_{1}(x)+{\bar{\psi}}_{2}(x))+{\bar{\psi}}_{1}(x)(1-{\bar{\psi}}_{2}(x))P(x).
\end{eqnarray*}
From the above definitions it follows that
\begin{eqnarray*}
{\bar{F}}_{n}(f,x)=\left\{\begin{array}{lr}
f(x),          &       x\in [0,x_1]\cup [x_4,1],   \\
f(x)(1-{\bar{\psi}}_{1}(x))+{\bar{\psi}}_{1}(x)
P(x),      &
x\in [x_1,x_2],  \\
P(x),          &       x\in [x_2,x_3],  \\
P(x)(1-{\bar{\psi}}_{2}(x))+{\bar{\psi}}_{2}(x)f(x), & x\in
[x_3,x_4].
            \end{array}
\right.
\end{eqnarray*}
Evidently, ${\bar{F}}_{n}$ is a positive linear operator which
depends on the functions values $f(k/n),\ 0 \leqslant k/n \leqslant
x_2$ or $x_3 \leqslant k/n \leqslant 1,$ it reproduces linear
functions, and ${\bar{F}}_{n} \in C^2([0,1])$ provided $f\in
W_\phi^{2}.$ Now for every $f\in C_{\bar{w}}$ define the Bernstein
type operator
\begin{eqnarray}
{\bar{B}}_{n}(f,x)&:=&B_{n}({\bar{F}}_{n}(f),x)  \nonumber\\
&=&\sum_{k/n \in [0,x_1]\cup [x_4,1]}p_{n,k}(x)f({\frac kn}) + \sum_{x_2<k/n<x_3}p_{n,k}(x)P({\frac kn})\nonumber\\
&&+\sum_{x_1<k/n<x_2}p_{n,k}(x)\{f({\frac kn})(1-{\bar{\psi}}_{1}({\frac kn}))+{\bar{\psi}}_{1}({\frac kn})P({\frac kn})\}\nonumber\\
&&+\sum_{x_3<k/n<x_4}p_{n,k}(x)\{P({\frac
kn})(1-{\bar{\psi}}_{2}({\frac kn}))+{\bar{\psi}}_{2}({\frac
kn})f({\frac kn})\}. \label{s1}
\end{eqnarray}
Obviously, ${\bar{B}}_{n}$ is a positive linear operator,
${\bar{B}}_{n}(f)$ is a polynomial of degree at most $n,$ it
preserves linear functions, and depends only on the function values
$f(k/n),\ k/n \in [0,x_2] \cup [x_3,1].$
Now we state our main results as follows:\\~\\
\textbf{Theorem 1.} {\it If \ $\alpha>0,$ for any $f \in
C_{\bar{w}},$ we have
\begin{eqnarray}
\|\bar{w}{\bar{B}}''_n(f)\|\leqslant Cn^2\|{\bar{w}}f\|.\label{s2}
\end{eqnarray}}
\textbf{Theorem 2.} {\it For any $\alpha >0,\
\min\{\beta(0),\beta(1)\} \geqslant {\frac 12},\ 0<\xi<1,$ we have
\begin{eqnarray}
|{\bar{w}(x)}\phi^2(x){\bar{B}}''_n(f,x)|\leqslant \left\{
\begin{array}{lrr}
Cn\|{\bar{w}}f\|,    &&f\in C_{\bar{w}},    \\
C\|{\bar{w}}\phi^{2}f''\|, &&f\in W_\phi^{2}.
              \end{array}\label{s3}
\right.
\end{eqnarray}}
\textbf{Theorem 3.} {\it For $f\in C_{\bar{w}},\ 0<\xi<1,\
\alpha>0,\ \min\{\beta(0),\beta(1)\} \geqslant {\frac 12},\ \alpha_0 \in (0,2),$ we have
\begin{eqnarray*}
{\bar{w}(x)}|f(x)-{\bar{B}_{n}(f,x)}|=O((n^{-{\frac
12}}\phi^{-1}(x)\delta_n(x))^{\alpha_0}) \Longleftrightarrow
\omega_\phi^2(f,t)_{\bar{w}}=O(t^{\alpha_0}).\label{s4}
\end{eqnarray*}}
%------------------------------------------------------------------------------------------
\section{Lemmas}
\textbf{Lemma 1.}(\cite{Zhou}) {\it For any non-negative real $u$ and
$v$, we have
\begin{eqnarray}
\sum_{k=1}^{n-1}({\frac kn})^{-u}(1-{\frac
kn})^{-v}p_{n,k}(x)\leqslant Cx^{-u}(1-x)^{-v}. \label{s5}
\end{eqnarray}}
\textbf{Lemma 2.}(\cite{Della}) {\it For any $\alpha >0,\ f \in
C_{\bar{w}},$ we have
\begin{eqnarray}
\|{\bar{w}}{\bar{B}_{n}(f)}\| \leqslant C\|{\bar{w}}f\|. \label{s6}
\end{eqnarray}}
%-----------------------------------------------------------------------
\textbf{Lemma 3.}(\cite{Wang}) {\it  Let $\min\{\beta(0), \beta(1)\}
\geqslant {\frac 12}$, then for $0<t<{\frac {1}{4}}$ and $t < x <
1-t,$ we have
\begin{eqnarray}
\int_{-{\frac {t}{2}}}^{\frac {t}{2}} \int_{-{\frac {t}{2}}}^{\frac
{t}{2}}\phi^{-2}(x+\sum_{k=1}^2u_k)du_1du_2 \leqslant
Ct^2\phi^{-2}(x).\label{s7}
\end{eqnarray}
} \textbf{Proof.} From the definition of $\phi(x),$ it is enough to
prove (\ref{s7}) for $t < x \leqslant {\frac 12}$ since the proof
for ${\frac 12} < x < 1-t$ is very similar. Obviously, we have
\begin{eqnarray*}
\int_{-{\frac {t}{2}}}^{\frac {t}{2}} \int_{-{\frac {t}{2}}}^{\frac
{t}{2}}{\frac {1}{(x+\sum_{k=1}^2u_k)}}du_1du_2 \leqslant
Ct^2x^{-1}.
\end{eqnarray*}
Therefore, by the H\"{o}lder inequality, we have
\begin{eqnarray*}
\int_{-{\frac {t}{2}}}^{\frac {t}{2}} \int_{-{\frac {t}{2}}}^{\frac
{t}{2}}\phi^{-2}(x+\sum_{k=1}^2u_k)du_1du_2\\
&\leqslant& C(\frac {16}{7})^{2\beta(1)}\int_{-{\frac
{t}{2}}}^{\frac {t}{2}} \int_{-{\frac {t}{2}}}^{\frac {t}{2}}{\frac
{1}{(x+\sum_{k=1}^2u_k)^{2\beta(0)}}}du_1du_2\\
&\leqslant& C(\frac
{16}{7})^{2\beta(1)}t^{2(1-2\beta(0))}(\int_{-{\frac {t}{2}}}^{\frac
{t}{2}} \int_{-{\frac {t}{2}}}^{\frac {t}{2}}{\frac
{1}{(x+\sum_{k=1}^2u_k)}}du_1du_2)^{2\beta(0)}\\
&\leqslant& C(\frac {16}{7})^{2\beta(1)}t^2x^{-2\beta(0)}. \Box
\end{eqnarray*}
%-------------------------------------------------------------------------------------------------------
\textbf{Lemma 4.}(\cite{Della}) {\it If $\gamma \in R,$ then
\begin{eqnarray}
\sum_{k=0}^np_{n,k}(x)|k-nx|^\gamma \leqslant Cn^{\frac
\gamma2}\varphi^\gamma(x).\label{s8}
\end{eqnarray}}
\textbf{Lemma 5.} {\it Let
$A_n(x):={\bar{w}(x)}\sum\limits_{|k-n\xi|\leqslant
\sqrt{n}}p_{n,k}(x)$. Then $A_n(x)\leqslant Cn^{-\frac \alpha2}$ for
$0<\xi <1$ and $\alpha>0$.} \\~\\
%----------------------------------------------------------------------
\textbf{Proof.} If $|x- \xi|\leqslant {\frac {3}{\sqrt{n}} }$, then
the statement is trivial. Hence assume $0 \leqslant x \leqslant
\xi-{\frac {3}{\sqrt{n}} }$ (the case $\xi+{\frac {3}{\sqrt{n}} }
\leqslant x \leqslant 1$ can be treated similarly). Then for a fixed
$x$ the maximum of $p_{n,k}(x)$ is attained for
$k=k_n:=[n\xi-\sqrt{n}]$. By using Stirling's formula, we get
\begin{eqnarray*}
p_{n,k_n}(x)&\leqslant& C{\frac {({\frac
{n}{e}})^n\sqrt{n}x^{k_n}(1-x)^{n-k_n}}{({\frac
{k_n}{e}})^{k_n}\sqrt{k_n}({\frac
{n-k_n}{e}})^{n-k_n}\sqrt{n-k_n}}}\nonumber\\
&\leqslant& {\frac {C}{\sqrt{n}}}({\frac {nx}{k_n}})^{k_n}({\frac
{n(1-x)}{n-k_n}})^{n-k_n}\nonumber\\
&=&{\frac {C}{\sqrt{n}}}(1-{\frac {k_n-nx}{k_n}})^{k_n}(1+{\frac
{k_n-nx}{n-k_n}})^{n-k_n}.
\end{eqnarray*}
Now from the inequalities
$$k_n-nx=[n\xi-\sqrt{n}]-nx>n(\xi-x)-\sqrt{n}-1\geqslant {\frac 12}n(\xi-x),$$
and
$$1-u\leqslant e^{-u-{\frac 12}u^2},\ 1+u\leqslant e^u,\ u\geqslant 0,$$
it follows that the second inequality is valid. To prove the first
one we consider the function $\lambda(u)=e^{-u-{\frac 12}u^2}+u-1.$
Here $\lambda(0)=0,\ \lambda^\prime(u)=-(1+u)e^{-u-{\frac
12}u^2}+1,\ \lambda^\prime(0)=0,\
\lambda^{\prime\prime}(u)=u(u+2)e^{-u-{\frac 12}u^2}\geqslant 0,$
whence $\lambda(u)\geqslant 0$ for $u\geqslant 0$. Hence
\begin{eqnarray*}
p_{n,k_n}(x) &\leqslant&{\frac {C}{\sqrt{n}}}exp\{k_n[-{\frac
{k_n-nx}{k_n}}-{\frac 12}({\frac {k_n-nx}{k_n}})^2] +
k_n-nx\}\nonumber\\
&=&{\frac {C}{\sqrt{n}}}exp\{-\frac {({k_n-nx})^2}{2k_n}\}\leqslant
e^{-Cn(\xi-x)^2}.
\end{eqnarray*}
Thus $A_n(x)\leqslant C(\xi-x)^\alpha e^{-Cn(\xi-x)^2}$. An easy
calculation shows that here the maximum is attained when
$\xi-x={\frac {C}{\sqrt{n}}}$ and the lemma follows. $\Box$\\~\\
\textbf{Lemma 6.} {\it For $0<\xi <1,\ \alpha,\ \beta>0$, we have
\begin{eqnarray}
{\bar{w}(x)}\sum\limits_{|k-n\xi|\leqslant \sqrt{n}}|k-nx|^\beta
p_{n,k}(x)\leqslant Cn^{\frac
{\beta-\alpha}{2}}\varphi^\beta(x).\label{s9}
\end{eqnarray}}
%-----------------------------------------------------------------------------------------------------
\textbf{Proof.} By (\ref{s8}) and the lemma 5, we have
\begin{eqnarray*}
{\bar{w}(x)}^{\frac
{1}{2n}}({\bar{w}(x)\sum\limits_{|k-n\xi|\leqslant
\sqrt{n}}p_{n,k}(x)})^{\frac
{2n-1}{2n}}(\sum\limits_{|k-n\xi|\leqslant \sqrt{n}}|k-nx|^{2n\beta}
p_{n,k}(x))^{\frac {1}{2n}}\leqslant Cn^{\frac
{\beta-\alpha}{2}}\varphi^\beta(x). \Box
\end{eqnarray*}
%------------------------------------------------------------------------------------------
\textbf{Lemma 7.} {\it For any $\alpha >0,\ f\in W_\phi^{2},\ \min\{\beta(0), \beta(1)\}
\geqslant {\frac 12},$ we
have
\begin{eqnarray}
{\bar{w}(x)}|f(x)-P(f,x)|_{[x_1,x_4]} \leqslant C({\frac
{\delta_n(x)}{\sqrt{n}\phi(x)}})^2\|{\bar{w}}\phi^{2}f''\|.\label{s10}
\end{eqnarray}}
\textbf{Proof.} If $x \in [x_1,x_4],$ for any $f\in W_\phi^{2},$ we
have
\begin{eqnarray*}
f(x_1)=f(x)+f'(x)(x_1-x)+\int_{x_1}^x (t-x_1)f''(t)dt,\\
f(x_4)=f(x)+f'(x)(x_4-x)+\int_{x_4}^x (t-x_4)f''(t)dt,\\
{\delta_n(x)} \sim {\frac {1}{\sqrt{n}}},\ n=1,2,\cdots.
\end{eqnarray*}
So
\begin{eqnarray*}
{\bar{w}(x)}|f(x)-P(f,x)| &\leqslant& {\bar{w}(x)}|{\frac
{x-x_4}{x_1-x_4}}|\int_{x_1}^x
|(t-x_1)f''(t)|dt\\
&+& {\bar{w}(x)}|{\frac
{x_1-x}{x_1-x_4}}|\int_{x_4}^x |(t-x_4)f''(t)|dt\\
&:=&I_1+I_2.
\end{eqnarray*}
Whence $t$ between $x_1$ and $x,$ we have
$\frac{|t-x_1|}{{\bar{w}}(t)}\leqslant
\frac{|x-x_1|}{{\bar{w}}(x)},$ then
\begin{eqnarray*}
I_1 \leqslant Cn^{\frac
12}\|{\bar{w}}\phi^{2}f''\||(x-x_1)(x-x_4)|\int_{x_1}^x\phi^{-2}(t)dt
&\leqslant&
C({\frac{\varphi(x)}{\sqrt{n}\phi(x)}})^2\|{\bar{w}}\phi^{2}f''\|\\
&\leqslant&
C({\frac{\delta_n(x)}{\sqrt{n}\phi(x)}})^2\|{\bar{w}}\phi^{2}f''\|.
\end{eqnarray*}
Analogously, we have
\begin{eqnarray*}
I_2 \leqslant
C({\frac{\delta_n(x)}{\sqrt{n}\phi(x)}})^2\|{\bar{w}}\phi^{2}f''\|.
\end{eqnarray*}
Now the lemma follows from combining these results together. $\Box$\\~\\
%------------------------------------------------------------------------------------------
\textbf{Lemma 8.} {\it If $f \in W_\phi^2,\ \min\{\beta(0),
\beta(1)\} \geqslant {\frac 12},$ then
\begin{eqnarray}
\|{\bar{w}}\phi^{2}{\bar{F}''_{n}}\|=O(\|{\bar{w}}\phi^{2}f''\|).\label{s11}
\end{eqnarray}}
\textbf{Proof.} Again, it is sufficient to estimate
$({\bar{w}}\phi^{2}{\bar{F}''_{n}})(x)$ for $x\in [x_3,x_4],$ and the
same as $x\in [x_1,x_2].$ For $x\in [x_2,x_3],$
${\bar{F}''_{n}}(x)=0,$ while for $x\in [0,x_1]\cup [x_4,1],$
${\bar{F}_{n}}(x)=f(x).$ Thus for $x\in [x_3,x_4],$ then
${\bar{F}_{n}}(x)=P(x)+{\bar{\psi}}_{2}(x)(f(x)-P(x))$ and
\begin{eqnarray*}
{\bar{F}''_{n}}(x)&=&n\psi''[n^{\frac
12}(x-x_3)](f(x)-P(x))\\
&& +\ 2n^{\frac 12}\psi'[n^{\frac
12}(x-x_3)](f(x)-P(x))'+\psi[n^{\frac 12}(x-x_3)]f''(x)\\
&& :=\ I_1(x) + I_2(x) + I_3(x).
\end{eqnarray*}
From the proof of lemma 7, we have
\begin{eqnarray*}
|{\bar{w}}(x)\phi^{2}(x)I_1(x)|&&\\
&=& O(n\phi^{2}(x)\psi''[n^{\frac
12}(x-x_3)]{\bar{w}}(x)(f(x)-P(x)))\\
&=& O(n\phi^{2}(x) \cdot
({\frac{\varphi(x)}{\sqrt{n}\phi(x)}})^2\|{\bar{w}}\phi^{2}f''\|)\\
&=& O(\|{\bar{w}}\phi^{2}f''\|).
\end{eqnarray*}
For $I_3(x),$ it is obvious that
\begin{eqnarray*}
|{\bar{w}}(x)\phi^{2}(x)I_3(x)|=O(\|{\bar{w}}\phi^{2}f''\|).
\end{eqnarray*}
Finally
\begin{eqnarray*}
|{\bar{w}}(x)\phi^{2}(x)I_2(x)|&&\\
&=& O(n^{\frac 12}{\bar{w}}(x)\phi^{2}(x)|f'(x)-P'(x)|)\\
&=& O(n^{\frac 12}{\bar{w}}(x)\phi^{2}(x)|f'(x)-n^{\frac
12}\int_{x_1}^{x_4}f'(t)dt|)\\
&=& O(n^{\frac 12}{\bar{w}}(x)\phi^{2}(x)|n^{\frac
12}\int_{x_1}^{x_4}\int_{t}^{x}f''(u)dudt|)\\
&=& O(n^{\frac
12}{\bar{w}}(x)\phi^{2}(x)|\int_{x_1}^{x_4}f''(u)du|)\\
&=& O(\|{\bar{w}}\phi^{2}f''\|). \Box
\end{eqnarray*}
%------------------------------------------------------------------------------------------
\section{Proof of Theorem}
\subsection{Proof of Theorem 1}
\textit{Case 1.}  If $f \in C_{\bar{w}},$ when $x\in [{\frac
1n},1-{\frac 1n}],$ by [2], we have
\begin{eqnarray}
|{\bar{w}(x)}{\bar{B}''_{n}(f,x)}| & \leqslant &
n\varphi^{-2}(x){\bar{w}(x)}|{\bar{B}_{n}(f,x)}|\nonumber\\
&+&{\bar{w}(x)}\varphi^{-4}(x)\sum_{k=0}^np_{n,k}(x)|k-nx||{\bar{F}}_{n}({\frac
kn})|\nonumber\\
& + &
{\bar{w}(x)}\varphi^{-4}(x)\sum_{k=0}^np_{n,k}(x)(k-nx)^2|{\bar{F}}_{n}({\frac
kn})|\nonumber\\  &:=& A_1+A_2+A_3.\label{s12}
\end{eqnarray}
By (\ref{s6}), we have
\begin{eqnarray}
A_1(x)=n\varphi^{-2}(x){\bar{w}(x)}|{\bar{B}_{n}(f,x)}| \leqslant
Cn^2\|{\bar{w}}f\|.\label{s13}
\end{eqnarray}
and
\begin{eqnarray*}
A_2={\bar{w}(x)}\varphi^{-4}(x)[\sum_{k/n\in
A}|k-nx||{\bar{F}}_{n}({\frac kn})|p_{n,k}(x)+\sum_{x_2 \leqslant
k/n\leqslant x_3}|k-nx||P({\frac
kn})|p_{n,k}(x)]:=\sigma_1+\sigma_2.
\end{eqnarray*}
thereof $A:=[0,x_2]\cup [x_3,1].$ If ${\frac kn}\in A,$ when ${\frac
{\bar{w}(x)}{\bar{w}(\frac {k}{n})}}\leqslant C(1+n^{-{\frac
{\alpha}{2}}}|k-nx|^\alpha),$ we have $|k-n\xi|\geqslant {\frac
{\sqrt{n}}{2}}$, by (\ref{s8}), then
\begin{eqnarray*}
\sigma_1 & \leqslant &
C\|{\bar{w}}f\|\varphi^{-4}(x)\sum_{k=0}^np_{n,k}(x)|k-nx|[1+n^{-{\frac
{\alpha}{2}}}|k-nx|^\alpha]\\
&=&
C\|{\bar{w}}f\|\varphi^{-4}(x)\sum_{k=0}^np_{n,k}(x)|k-nx|+Cn^{-{\frac
{\alpha}{2}}}\|{\bar{w}}f\|\varphi^{-4}(x)\sum_{k=0}^np_{n,k}(x)|k-nx|^{1+\alpha}\\
& \leqslant & Cn^{\frac 12}\varphi^{-3}(x)\|{\bar{w}}f\|+Cn^{\frac
12}\varphi^{-3+\alpha}(x)\|{\bar{w}}f\| \\
& \leqslant & Cn^2\|{\bar{w}}f\|.
\end{eqnarray*}
For $\sigma_2,$ $P$ is a linear function. We note $|P({\frac
kn})|\leqslant max(|P(x_1)|,|P(x_4)|):=P(a).$ If $x\in [x_1,x_4],$
we have ${\bar{w}(x)}\leqslant {\bar{w}(a)}.$ So, if $x\in
[x_1,x_4],$ by (\ref{s8}), then
\begin{eqnarray*}
\sigma_2\leqslant
C{\bar{w}(a)}P(a)\varphi^{-4}(x)\sum_{k=0}^np_{n,k}(x)|k-nx|\leqslant
Cn^{2}\|{\bar{w}}f\|.
\end{eqnarray*}
%----------------------------------------------------------------------------------------------------------------
If $x\notin [x_1,x_4],$ then ${\bar{w}(a)}>n^{-{\frac
{\alpha}{2}}},$ by (\ref{s9}), we have
\begin{eqnarray*}
\sigma_2 &\leqslant& C{\bar{w}(x)}\varphi^{-4}(x) \sum_{x_2
\leqslant k/n\leqslant x_3}|P(a)||(k-nx)|p_{n,k}(x)\\
&\leqslant& Cn^{\frac
{\alpha}{2}}\|{\bar{w}}f\|\varphi^{-4}(x){\bar{w}(x)}\sum_{x_2
\leqslant k/n\leqslant x_3}|k-nx|p_{n,k}(x)\\
&\leqslant& Cn^2\|{\bar{w}}f\|.
\end{eqnarray*}
So
\begin{eqnarray}
A_2 \leqslant Cn^2\|{\bar{w}}f\|.\label{s14}
\end{eqnarray}
Similarly
\begin{eqnarray}
A_3 \leqslant Cn^2\|{\bar{w}}f\|.\label{s15}
\end{eqnarray}
It follows from combining with (\ref{s12})-(\ref{s15}) that the inequality is proved.\\
\textit{Case 2.} When $x\in [0,{\frac 1n}]$ (The same as $x\in
[1-{\frac 1n},1]$), by \cite{Totik}, then
\begin{eqnarray*}
{\bar{B}''_{n}(f,x)}=n(n-1)\sum_{k=0}^{n-2}\overrightarrow{\Delta}_{\frac
1{n}}^{2}{\bar{F}}_{n}{(\frac kn)}p_{n-2,k}(x).
\end{eqnarray*}
We have
\begin{eqnarray*}
|{\bar{w}(x)}{\bar{B}''_{n}(f,x)}| &\leqslant&
Cn^2{\bar{w}(x)}\sum_{k=0}^{n-2}|\overrightarrow{\Delta}_{\frac
1{n}}^{2}{\bar{F}}_{n}{(\frac kn)}|p_{n-2,k}(x)\\
&=& Cn^2{\bar{w}(x)}[\sum_{k/n\in
A}p_{n-2,k}(x)|\overrightarrow{\Delta}_{\frac
1{n}}^{2}{\bar{F}}_{n}{(\frac kn)}|+\sum_{x_2 \leqslant k/n\leqslant
x_3}p_{n-2,k}(x)|\overrightarrow{\Delta}_{\frac 1{n}}^{2}P({\frac
kn})|].
\end{eqnarray*}
We can deal with it in accordance with Case 1,
and prove it immediately, then the theorem is done. $\Box$\\

%------------------------------------------------------------------------------------------
\subsection{Proof of Theorem 2}
 (1) We prove the first inequality of Theorem 2.\\
 \textit{Case 1.} If $0\leqslant \varphi(x)\leqslant {\frac
{1}{\sqrt{n}}}$, by (\ref{s2}), we have
\begin{eqnarray*}
|{\bar{w}(x)}\phi^{2}(x){\bar{B}}''_n(f,x)|= \varphi^2(x)\cdot{\frac
{\phi^{2}(x)}{\varphi^2(x)}}|{\bar{w}(x)}{\bar{B}}''_n(f,x)|\leqslant
Cn\|{\bar{w}}f\|.
\end{eqnarray*}
\textit{Case 2.} If $\varphi(x)> {\frac {1}{\sqrt{n}}}$, by
\cite{Totik},\ we have
\begin{eqnarray*}
{\bar{B}}''_n(f,x)=B''_n({\bar{F}_{n}},x) & = &
(\varphi^{2}(x))^{-1}\sum_{i=0}^{2}Q_{i}(x,n)n^i\sum_{k=0}^{n}(x-{\frac
kn})^{i}{\bar{F}}_{n}({\frac kn})p_{n,k}(x),\\
(\varphi^{2}(x))^{-1}Q_{i}(x,n)n^{i} &\leqslant&
C(n/\varphi^{2}(x))^{1+i/2}.
\end{eqnarray*}
So
\begin{eqnarray*}
&&|{\bar{w}}(x)\phi^{2}(x){\bar{B}}''_n(f,x)|\nonumber\\
&\leqslant& C{\bar{w}(x)}\phi^{2}(x)\sum_{i=0}^{2}({\frac
{n}{\varphi^2(x)}})^{1+i/2}\sum_{k=0}^{n}|x-{\frac
kn}|^{i}|{\bar{F}}_{n}({\frac kn})|p_{n,k}(x)\nonumber\\
&=& C{\bar{w}(x)}\phi^{2}(x)\sum_{i=0}^{2}({\frac
{n}{\varphi^2(x)}})^{1+i/2}\sum_{k/n\in A}|x-{\frac
kn}|^{i}|{\bar{F}}_{n}({\frac kn})|p_{n,k}(x)\nonumber\\
&+& C{\bar{w}(x)}\phi^{2}(x)\sum_{i=0}^{2}({\frac
{n}{\varphi^2(x)}})^{1+i/2}\sum_{x_2 \leqslant k/n \leqslant
x_3}|{x-{\frac kn}}|^{i}|P({\frac
kn})|p_{n,k}(x)\nonumber\\
&:=&\sigma_1+ \sigma_2.
\end{eqnarray*}
Where $A:=[0,x_2]\cup [x_3,1]$. Working as in the proof of Theorem
1, We can get $\sigma_1\leqslant Cn\|{\bar{w}}f\|,$
$\sigma_2\leqslant Cn\|{\bar{w}}f\|.$ By bringing these
facts together, we can immediately get the first inequality of Theorem 2.\\
(2) If $f\in W_\phi^{2},$ by (\ref{s1}), then
\begin{eqnarray}
|{\bar{w}(x)}\phi^{2}(x){\bar{B}}''_n(f,x)| & \leqslant &
n^2{\bar{w}(x)}\phi^{2}(x)\sum_{k=0}^{n-2}|\overrightarrow{\Delta}_{\frac
1{n}}^{2}{\bar{F}}_{n}{(\frac kn)}|p_{n-2,k}(x)\nonumber\\
&=&
n^2{\bar{w}(x)}\phi^{2}(x)\sum_{k=1}^{n-3}|\overrightarrow{\Delta}_{\frac
1{n}}^{2}{\bar{F}}_{n}{(\frac kn)}|p_{n-2,k}(x)\nonumber\\
&+& n^2{\bar{w}(x)}\phi^{2}(x)|\overrightarrow{\Delta}_{\frac
1{n}}^{2}{\bar{F}}_{n}{(0)}|p_{n-2,0}(x)\nonumber\\
&+& n^2{\bar{w}(x)}\phi^{2}(x)|\overrightarrow{\Delta}_{\frac
1{n}}^{2}{\bar{F}}_{n}{(\frac {n-2}{n})}|p_{n-2,n-2}(x) \nonumber\\
&:=& I_1+I_2+I_3.\label{s16}
\end{eqnarray}
By \cite{Totik}, if $0<k<n-2,$ we have
\begin{eqnarray}
|\overrightarrow{\Delta}_{\frac 1{n}}^{2}{\bar{F}}_{n}({\frac kn})
|\leqslant Cn^{-1}\int_{0}^{\frac {2}{n}}|{\bar{F}}''_{n}({\frac
kn}+u)|du,\label{s17}
\end{eqnarray}
%---------------------------------------------------------------------
If $k=0,$ we have
\begin{eqnarray}
|\overrightarrow{\Delta}_{\frac 1{n}}^{2}{\bar{F}}_{n}(0)| \leqslant
C\int_{0}^{\frac {2}{n}}u|{\bar{F}}''_{n}(u)|du,\label{s18}
\end{eqnarray}
Similarly
\begin{eqnarray}
|\overrightarrow{\Delta}_{\frac 1n}^2{\bar{F}}_{n}({\frac
{n-2}{n}})| &\leqslant& Cn^{-1}\int_{1-{\frac
{2}{n}}}^{1}(1-u)|{\bar{F}}''_{n}(u)|du.\label{s19}
\end{eqnarray}
By (\ref{s17}), then
\begin{eqnarray}
I_1 &\leqslant&
Cn{\bar{w}(x)}\phi^{2}(x)\sum_{k=1}^{n-3}\int_{0}^{\frac
{2}{n}}|{\bar{F}}''_{n}({\frac kn}+u)|dup_{n-2,k}(x) \nonumber\\
&=& Cn{\bar{w}(x)}\phi^{2}(x)\sum_{k/n\in A}\int_{0}^{\frac
{2}{n}}|{\bar{F}}''_{n}({\frac
kn}+u)|dup_{n-2,k}(x) \nonumber\\
&+& Cn{\bar{w}(x)}\phi^{2}(x)\sum_{x_2 \leqslant k/n \leqslant
x_3}\int_{0}^{\frac {2}{n}}|P''({\frac
kn}+u)|dup_{n-2,k}(x) \nonumber\\
&:=& T_1+T_2.\nonumber
\end{eqnarray}
Where $A:=[0,x_2]\cup [x_3,1],$ $P$ is a linear function. If ${\frac
kn}\in A,$ when ${\frac {\bar{w}(x)}{\bar{w}(\frac
{k}{n})}}\leqslant C(1+n^{-{\frac {\alpha}{2}}}|k-nx|^\alpha),$ we
have $|k-n\xi|\geqslant {\frac {\sqrt{n}}{2}}$, by (\ref{s5}),
(\ref{s8}) and (\ref{s11}), then
\begin{eqnarray*}
T_1 &\leqslant&
C{\bar{w}(x)}\phi^{2}(x)\|{\bar{w}}\phi^{2}{\bar{F}''_{n}}\|\sum_{k/n\in
A}p_{n-2,k}(x){\bar{w}}^{-1}({\frac kn})\phi^{-2}({\frac kn})\\
&\leqslant&
C\phi^{2}(x)\|{\bar{w}}\phi^{2}{\bar{F}''_{n}}\|\sum_{k=0}^{n-2}p_{n-2,k}(x)(1+n^{-{\frac
{\alpha}{2}}}|k-nx|^\alpha)\phi^{-2}({\frac kn})\\
&\leqslant& C\|{\bar{w}}\phi^{2}{\bar{F}''_{n}}\|\\
&\leqslant& C\|{\bar{w}}\phi^{2}f''\|.
\end{eqnarray*}
Working as the Theorem 1, we can get
\begin{eqnarray*}
T_2 \leqslant C\|{\bar{w}}\phi^{2}f''\|.
\end{eqnarray*}
So, we can get
\begin{eqnarray}
I_1 \leqslant C\|{\bar{w}}\phi^{2}f''\|.\label{s20}
\end{eqnarray}
By (\ref{s11}) and (\ref{s18}), we have
\begin{eqnarray}
I_2 &\leqslant&
Cn^2{\bar{w}(x)}\phi^{2}(x)(1-x)^{n-2}\int_{0}^{\frac
{2}{n}}u|{\bar{F}}''_{n}(u)|du  \nonumber\\
&\leqslant&
Cn^2{\bar{w}(x)}\phi^{2}(x)(1-x)^{n-2}\|{\bar{w}}\phi^{2}{\bar{F}''_{n}}\|\int_{0}^{\frac
{2}{n}}u{\bar{w}}^{-1}(u)\phi^{-2}(u)du  \nonumber\\
&\leqslant& C\|{\bar{w}}\phi^{2}{\bar{F}''_{n}}\|  \nonumber\\
&\leqslant& C\|{\bar{w}}\phi^{2}f''\|.  \label{s21}
\end{eqnarray}
Similarly
\begin{eqnarray}
I_3 \leqslant C\|{\bar{w}}\phi^{2}f''\|.\label{s22}
\end{eqnarray}
By bringing (\ref{s16}), (\ref{s20})-(\ref{s22}) together, we
can get the second inequality
of Theorem 2. $\Box$\\~\\
%-------------------------------------------------------------------------------------------------------------
 \textbf{Corollary} {\it For any $\alpha >0,\ 0 \leqslant \lambda
\leqslant 1,$ we have
\begin{eqnarray}
|{\bar{w}(x)}\varphi^{2\lambda}(x){\bar{B}}''_n(f,x)|\leqslant
\left\{
\begin{array}{lrr}
Cn{\{max\{n^{1-\lambda},\varphi^{2(\lambda-1)}\}\}}\|{\bar{w}}f\|,    &&f\in C_{\bar{w}},    \\
C\|{\bar{w}}\varphi^{2\lambda}f''\|, &&f\in
W_{{\bar{w}},\lambda}^{2}.
              \end{array}\label{s23}
\right.
\end{eqnarray}}

%---------------------------------------------------------------------------------------------------------------------------
\subsection{Proof of Theorem 3}
\subsubsection{The direct theorem}
We know
\begin{eqnarray}
{\bar{F}}_n(t)={\bar{F}}_n(x)+{\bar{F}}'_n(t)(t-x)+\int_x^t (t-u){\bar{F}}''_n(u)du, \label{s24}\\
B_n(t-x,x)=0.\label{s25}
\end{eqnarray}
According to the definition of $W_\phi^{2},$ by (\ref{s24}) and (\ref{s25}), for any $g \in
W_\phi^{2},$ we have
${\bar{B}}_{n}(g,x)=B_{n}({\bar{G}}_{n}(g),x),$ then
\begin{eqnarray}
{\bar{w}(x)}|{\bar{G}}_{n}(x)-B_{n}({\bar{G}}_{n},x)|={\bar{w}(x)}|B_n(R_2({\bar{G}}_n,t,x),x)|,\label{s26}
\end{eqnarray}
thereof $R_2({\bar{G}}_n,t,x)=\int_x^t (t-u){\bar{G}}''_n(u)du.$
\begin{eqnarray}
{\bar{w}(x)}|{\bar{G}}_{n}(x)-B_{n}({\bar{G}}_{n},x)| &\leqslant&
C{\bar{w}(x)}\sum_{k=1}^{n-1}p_{n,k}(x)\int_x^{\frac kn}|{\frac
kn}-u||{\bar{G}}''_n(u)|du\nonumber\\
&+& C{\bar{w}(x)}p_{n,0}(x)\int_0^xu|{\bar{G}}''_n(u)|du\nonumber\\
&+&
C{\bar{w}(x)}p_{n,n}(x)\int_x^1(1-u)|{\bar{G}}''_n(u)|du\nonumber\\
&:=& I_1+I_2+I_3.\label{s27}
\end{eqnarray}
If $u$ between ${\frac kn}$ and $x,$ we have
\begin{eqnarray}
\frac{|{\frac kn}-u|}{{\bar{w}}^2(u)}\leqslant \frac{|{\frac
kn}-x|}{{\bar{w}}^2(x)},\ \frac{|{\frac
kn}-u|}{\phi^{4}(u)}\leqslant
\frac{|{\frac
kn}-x|}{\phi^{4}(x)}.\label{s28}
\end{eqnarray}
By (\ref{s8}) and (\ref{s28}), then
\begin{eqnarray}
I_1 &\leqslant&
C\|{\bar{w}}\phi^{2}{\bar{G}}''_n\|{\bar{w}(x)}\sum_{k=1}^{n-1}p_{n,k}(x)\int_x^{\frac
kn}{\frac {|{\frac
kn}-u|}{{\bar{w}(u)}\phi^{2}(u)}}du  \nonumber\\
&\leqslant&
C\|{\bar{w}}\phi^{2}{\bar{G}}''_n\|{\bar{w}(x)}\sum_{k=1}^{n-1}p_{n,k}(x)(\int_x^{\frac
kn}{\frac {|{\frac
kn}-u|}{\bar{w}^2(u)}}du)^{\frac 12}({\int_x^{\frac
kn}{\frac
{|{\frac
kn}-u|}{\phi^{4}(u)}}du})^{\frac 12}  \nonumber\\
&\leqslant&Cn^{-2}\|{\bar{w}}\phi^{2}{\bar{G}}''_n\|\phi^{-2}(x)\sum_{k=0}^{n-1}p_{n,k}(x)(k-nx)^2  \nonumber\\
&\leqslant& Cn^{-1}{\frac
{\varphi^{2}(x)}{\phi^{2}(x)}}\|{\bar{w}}\phi^{2}{\bar{G}}''_n\| \nonumber\\
&\leqslant& Cn^{-1}{\frac
{\delta_n^{2}(x)}{\phi^{2}(x)}}\|{\bar{w}}\phi^{2}{\bar{G}}''_n\|  \nonumber\\
&=& C(\frac
{\delta_n(x)}{\sqrt{n}\phi(x)})^2\|{\bar{w}}\phi^{2}{\bar{G}}''_n\|.
\label{s29}
\end{eqnarray}
For $I_2,$ when $u$ between ${\frac kn}$ and $x,$ we let $k=0,$ then
$\frac{u}{{\bar{w}}(u)}\leqslant \frac{x}{{\bar{w}}(x)},$ and
\begin{eqnarray}
I_2 &\leqslant&
C\|{\bar{w}}\phi^{2}{\bar{G}}''_n\|{\bar{w}(x)}p_{n,0}(x)\int_0^xu{\bar{w}^{-1}(u)}\phi^{-2}(u)du\nonumber\\
&\leqslant& C(nx)(1-x)^{n-1} \cdot n^{-1}{\frac
{\varphi^{2}(x)}{\phi^{2}(x)}}\|{\bar{w}}\phi^{2}{\bar{G}}''_n\|\nonumber\\
&\leqslant& C(\frac
{\delta_n(x)}{\sqrt{n}\phi(x)})^2\|{\bar{w}}\phi^{2}{\bar{G}}''_n\|.\label{s30}
\end{eqnarray}
Similarly, we have
\begin{eqnarray}
I_3 \leqslant C(\frac
{\delta_n(x)}{\sqrt{n}\phi(x)})^2\|{\bar{w}}\phi^{2}{\bar{G}}''_n\|.\label{s31}
\end{eqnarray}
By bringing (\ref{s29})-(\ref{s31}), we have
\begin{eqnarray}
{\bar{w}(x)}|{\bar{G}}_{n}(x)-B_{n}({\bar{G}}_{n},x)| \leqslant
C(\frac
{\delta_n(x)}{\sqrt{n}\phi(x)})^2\|{\bar{w}}\phi^{2}{\bar{G}}''_n\|.\label{s32}
\end{eqnarray}
By (\ref{s10}) and (\ref{s32}), when $g \in
W_\phi^{2},$ then
\begin{eqnarray}
{\bar{w}(x)}|g(x)-{\bar{B}_{n}(g,x)}| &\leqslant&
{\bar{w}(x)}|g(x)-{\bar{G}}_{n}(g,x)| +
{\bar{w}(x)}|{\bar{G}}_{n}(g,x)-{\bar{B}_{n}(g,x)}|\nonumber\\
&\leqslant& {\bar{w}(x)}|g(x)-P(g,x)|_{[x_1,x_4]} + C(\frac
{\delta_n(x)}{\sqrt{n}\phi(x)})^2\|{\bar{w}}\phi^{2}{\bar{G}}''_n\|\nonumber\\
&\leqslant& C(\frac
{\delta_n(x)}{\sqrt{n}\phi(x)})^2\|{\bar{w}}\phi^{2}g''\|.\label{s33}
\end{eqnarray}
For $f \in C_{\bar{w}},$ we choose proper $g \in W_\phi^{2},$ by
(\ref{s6}) and (\ref{s33}), then
\begin{eqnarray*}
{\bar{w}(x)}|f(x)-{\bar{B}_{n}(f,x)}| &\leqslant&
{\bar{w}(x)}|f(x)-g(x)| + {\bar{w}(x)}|{\bar{B}_{n}(f-g,x)}| +
{\bar{w}(x)}|g(x)-{\bar{B}_{n}(g,x)}|\\
&\leqslant& C(\|{\bar{w}}(f-g)\|+(\frac
{\delta_n(x)}{\sqrt{n}\phi(x)})^2\|{\bar{w}}\phi^{2}g''\|)\\
&\leqslant& C\omega_\phi^2(f,\frac
{\delta_n(x)}{\sqrt{n}\phi(x)})_{\bar{w}}. \Box
\end{eqnarray*}

\subsubsection{The inverse theorem}
The main-part K-functional is
given by
\begin{eqnarray*}
K_{2,\phi}(f,t^2)_{\bar{w}}=\sup_{0<h \leqslant
t}\inf_g\{\|{\bar{w}}(f-g)\|+t^2\|{\bar{w}}\phi^{2}g''\|,\ g' \in
A.C._{loc}\}.
\end{eqnarray*}
By \cite{Totik}, we have
\begin{eqnarray}
C^{-1}K_{2,\phi}(f,t^2)_{\bar{w}} \leqslant
\omega_\phi^2(f,t)_{\bar{w}} \leqslant
CK_{2,\phi}(f,t^2)_{\bar{w}}.\label{s34}
\end{eqnarray}
\textbf{Proof.} Let $\delta>0,$ by (\ref{s34}), we choose proper $g$
so that
\begin{eqnarray}
\|{\bar{w}}(f-g)\| \leqslant C\omega_\phi^2(f,\delta)_{\bar{w}},\
\|{\bar{w}}\phi^{2}g''\| \leqslant
C\delta^{-2}\omega_\phi^2(f,\delta)_{\bar{w}}.\label{s35}
\end{eqnarray}
then
\begin{eqnarray}
|{\bar{w}}(x)\Delta_{h\phi}^2f(x)| &\leqslant&
|{\bar{w}}(x)\Delta_{h\phi}^2(f(x)-{\bar{B}_{n}(f,x)})|+|{\bar{w}}(x)\Delta_{h\phi}^2\bar{B}_{n}(f-g,x)|\nonumber\\
&+& |{\bar{w}}(x)\Delta_{h\phi}^2{\bar{B}_{n}(g,x)}|\nonumber\\
&\leqslant& \sum_{j=0}^2C_2^j(n^{-\frac
12}{\frac {\delta_n(x+(1-j)h\phi(x))}{{\phi(x+(1-j)h\phi(x))}}})^{\alpha_0}\nonumber\\
&+& \int_{-{\frac {h\phi(x)}{2}}}^{\frac {h\phi(x)}{2}}
\int_{-{\frac {h\phi(x)}{2}}}^{\frac
{h\phi(x)}{2}}{\bar{w}}(x){\bar{B}''_{n}(f-g,x+\sum_{k=1}^2u_k)}du_1
du_2\nonumber\\
&+& \int_{-{\frac {h\phi(x)}{2}}}^{\frac {h\phi(x)}{2}}
\int_{-{\frac {h\phi(x)}{2}}}^{\frac
{h\phi(x)}{2}}{\bar{w}}(x){\bar{B}''_{n}(g,x+\sum_{k=1}^2u_k)}du_1
du_2\nonumber\\
&:=& J_1+J_2+J_3.\label{s36}
\end{eqnarray}
Obviously
\begin{eqnarray}
J_1 \leqslant C((n^{-{\frac
12}}\phi^{-1}(x)\delta_n(x))^{\alpha_0}).\label{s37}
\end{eqnarray}
By (\ref{s2}) and (\ref{s35}), we have
\begin{eqnarray}
J_2 &\leqslant& Cn^2\|{\bar{w}}(f-g)\|\int_{-{\frac
{h\phi(x)}{2}}}^{\frac {h\phi(x)}{2}} \int_{-{\frac
{h\phi(x)}{2}}}^{\frac
{h\phi(x)}{2}}du_1du_2\nonumber\\
&\leqslant& Cn^2h^2\phi^{2}(x)\|{\bar{w}}(f-g)\|\nonumber\\
&\leqslant&
Cn^2h^2\phi^{2}(x)\omega_\phi^2(f,\delta)_{\bar{w}}.\label{s38}
\end{eqnarray}
By the second inequality of (\ref{s23}) and (\ref{s35}),\ we have
\begin{eqnarray}
J_2 &\leqslant& Cn\|{\bar{w}}(f-g)\|\int_{-{\frac
{h\phi(x)}{2}}}^{\frac {h\phi(x)}{2}} \int_{-{\frac
{h\phi(x)}{2}}}^{\frac
{h\phi(x)}{2}}\varphi^{-2}(x+\sum_{k=1}^2u_k)du_1du_2\nonumber\\
&\leqslant& Cnh^2\phi^2(x)\varphi^{-2}(x)\|{\bar{w}}(f-g)\|\nonumber\\
&\leqslant&
Cnh^2\phi^2(x)\varphi^{-2}(x)\omega_\phi^2(f,\delta)_{\bar{w}}.\label{s39}
\end{eqnarray}
By the second inequality of (\ref{s3}), (\ref{s7})and (\ref{s35}),\ we have
\begin{eqnarray}
J_3 &\leqslant& C\|{\bar{w}}\phi^{2}g''\|{\bar{w}(x)}\int_{-{\frac
{h\phi(x)}{2}}}^{\frac {h\phi(x)}{2}} \int_{-{\frac
{h\phi(x)}{2}}}^{\frac
{h\phi(x)}{2}}{\bar{w}^{-1}(x+\sum_{k=1}^2u_k)}\phi^{-2}(x+\sum_{k=1}^2u_k)du_1du_2\nonumber\\
&\leqslant& Ch^2\|{\bar{w}}\phi^{2}g''\|\nonumber\\
&\leqslant&
Ch^2\delta^{-2}\omega_\phi^2(f,\delta)_{\bar{w}}.\label{s40}
\end{eqnarray}
Now, by (\ref{s36})-(\ref{s40}), there
exists a constant $M>0$ so that
\begin{align*}
|\bar{w}(x)\Delta_{h\phi}^2f(x)| \leqslant C((n^{-\frac
12}{\frac {\delta_n(x)}{\phi(x)}})^{\alpha_0}\\
+ \min\{n{\frac
{\phi^2(x)}{\varphi^2(x)}},n^2\phi^2(x)\}h^2\omega_{\phi}^{2}(f,\delta)_{\bar{w}}
+ h^2\delta^{-2}\omega_{\phi}^{2}(f,\delta)_{\bar{w}})\\
\leqslant C((n^{-\frac
12}{\frac {\delta_n(x)}{\phi(x)}})^{\alpha_0}\\
+\  h^2M^2(n^{-\frac 12}{\frac {\varphi(x)}{\phi(x)}}\ +\ n^{-\frac
12}{\frac
{n^{-1/2}}{\phi(x)}})^{-2}\omega_{\phi}^{2}(f,\delta)_{\bar{w}}
+ h^2\delta^{-2}\omega_{\phi}^{2}(f,\delta)_{\bar{w}})\\
\leqslant C((n^{-\frac
12}{\frac {\delta_n(x)}{\phi(x)}})^{\alpha_0}\\
+ h^2M^2(n^{-\frac 12}{\frac
{\delta_n(x)}{\phi(x)}})^{-2}\omega_{\phi}^{2}(f,\delta)_{\bar{w}}
+ h^2\delta^{-2}\omega_{\phi}^{2}(f,\delta)_{\bar{w}}).\\
\end{align*}
When $n \geqslant 2,$ we have
\begin{eqnarray*}
n^{-\frac 12}\delta_n(x) < (n-1)^{-\frac 12}\delta_{n-1}(x)
\leqslant \sqrt{2}n^{-\frac 12}\delta_n(x),
\end{eqnarray*}
Choosing proper $x,\ \delta,\  n \in N,$ so that
\begin{align*}
n^{-\frac 12}{\frac {\delta_n(x)}{\phi(x)}} \leqslant \delta <
(n-1)^{-\frac 12}{\frac {\delta_{n-1}(x)}{\phi(x)}},
\end{align*}
Therefore
\begin{align*}
|{\bar{w}}(x)\Delta_{h\phi}^2f(x)| \leqslant C\{\delta^{\alpha_0} +
h^2\delta^{-2}\omega_{\phi}^{2}(f,\delta)_{\bar{w}}\}.
\end{align*}
Which implies
\begin{align*}
\omega_\phi^{2}(f,t)_{\bar{w}} \leqslant C\{\delta^{\alpha_0} +
h^2\delta^{-2}\omega_{\phi}^{2}(f,\delta)_{\bar{w}}\}.\label{s38}
\end{align*}
So, by Berens-Lorentz lemma in \cite{Totik}, we get
\begin{align*}
\omega_\phi^{2}(f,t)_{\bar{w}} \leqslant Ct^{\alpha_0}. \Box
\end{align*}

%---------------------------------------------------------------------------------------------------------------------------------------------

\end{document}